\documentclass{article}
\usepackage[T2A]{fontenc}
\usepackage[utf8]{inputenc}
\usepackage[english,russian]{babel}
\usepackage[tbtags]{amsmath}
\usepackage{amsfonts,amssymb,mathrsfs,amscd}
\usepackage[hyper]{msb-a}

\makeatletter
\gdef\No{{\select@language{russian}\textnumero}}
\makeatother

\JournalName{}
\numberwithin{equation}{section}
\theoremstyle{plain}

\theoremstyle{definition}

\def \Z {{\mathbf {Z}}}
\def \J {{\mathbf {J}}}
\def \N {{\mathbf {N}}}

\def \B {{\cal B}}
\def \T {{\mathbf  T}}
\def\uu{\bigsqcup}
\def\eps{\varepsilon}

\begin{document}


\title{Типичные расширения эргодических систем}
\author[V.\,V.~Ryzhikov]{В.\,В.~Рыжиков}
\address{Московский государственный университет}
\email{vryzh@mail.ru}

\date{10.10.2022}
\udk{517.987}

\maketitle

\begin{fulltext}

\begin{abstract}
Статья посвящена задачам  о типичных свойствах расширений динамических систем  с инвариантной  мерой.  Доказано, что типичные расширения сохраняют сингулярность спектра, свойство перемешивания и некоторые другие асимптотические свойства. Обнаружено, что   сохранение алгебраических свойств, вообще говоря, зависит от статистических свойств базы. Установлено, что $P$-энтропия типичного расширения   принимает бесконечное значение. Это дает новое доказательство  результа Вейса, Глазнера, Остина, Тувено о недоминантности  детерминированных действий.   Рассмотрены   типичные измеримые семейства автоморфизмов вероятностного пространства. В асимптотическом  поведении  представителей типичного  семейства  показан  их динамический конформизм вместе  динамическим индивидуализмом.

Библиография: 15 названий.
\end{abstract}


\markright{Типичные расширения эргодических систем}


\section{Введение} 
Группа автоморфизмов $ Aut=Aut(\mu)$ стандартного вероятностного пространства $(X,\B,\mu)$ оснащается  полной метрикой Халмоша $\rho $. Расстояние между автоморфизмами $S$ и $T$ определяется формулой  
$$ \rho(S,T)=\sum_i 2^{-i}\left(\mu(SA_i\Delta TA_i)+\mu(S^{-1}A_i\Delta T^{-1}A_i)\right),$$
где $\{A_i\}$  -- некоторое  фиксированное семейство  множеств,  плотное в алгебре $\B$.
Говорят, что множество автоморфизмов типично, если оно содержит  некоторое $G_\delta$-множество, плотное в  $Aut$. 
Говорят, что свойство автоморфизма типично, если множество  автоморфизмов, обладающих этим свойством, типично. Когда специалисты традиционно говорят, что "типичный автоморфизм" является эргодическим, на самом деле они подразумевают в точности то, что множество эргодических автоморфизмов является типичным. Мы также будем использовать подобный математический жаргон.

Теория типичных действий с инвариантной мерой имеет давнюю историю,
находит приложения, например,  для спектральной теории динамических систем (см.\cite{ST}) и по-прежнему привлекает внимание исследователей, о чем свидетельствуют недавние работы \cite{Bo}--\cite{AGTW}.В статье \cite{GTW} показано, что типичные расширения сохраняют $K$-свойство и бернуллиевость базы,  в   работе \cite{AGTW} установлен изоморфизм   эргодического преобразования $S$ с положительной энтропией своему типичному расширению (слово \it преобразование \rm мы используем как  синоним термина \it автоморфизм\rm). Преобразование  с нулевой энтропией, напротив,  не изоморфно своему  типичному расширению  \cite{AGTW}.  Мы дадим другое доказательство этого факта, 
используя  метод работы \cite{R21}: для заданного автоморфизма $S$ подбираем числовой инвариант $h_P$ (типа энтропии Кушниренко \cite{Ku}) такой, что $h_P(S)=0$, но  для  типичных  расширений $R$ автоморфизма $S$ выполнено $h_P(R)=\infty$. 

Автоморфизм,  сопряженный своему обратному, называется  симметричным.
Свойство симметричности сохраняется для действий с положительной энтропией, что сразу вытекает из результата \cite{AGTW} об  изоморфизме. Мы покажем, что это свойство  не сохраняется при типичных расширениях жесткого симметричного автоморфизма.  Некоторые спектральные и асимптотические инварианты наследуются  типичными  расширениями:  мы покажем, что таковыми свойствами являются   сингулярность    спектра автоморфизма,  
частичная  жесткость, мягкое и сильное перемешивание. 

Теорию типичных действий формально можно рассматривать как  теорию типичных расширений тождественного действия на одноточечном пространстве. Если базовое пространство с мерой состоит из конечного набора точек, то  при рассмотрении  расширений тождественного действия,   мы изучаем   типичные свойства  конечного набора действий. Типичные расширения тождественного преобразования, заданного  на стандартном пространстве Лебега, характеризуют типичные  измеримые континуальные семейства преобразований. В частности, будет показано следующее. Для типичного  семейства $\{T_x\}$ найдется  последовательность  $m_i\to\infty$ такая, что степени  
$T^{m_i}_x$ сходятся к тождественному оператору (динамический конформизм). При этом для некоторой другой последовательности $n_j\to\infty$ степени $T^{n_i}_x$ слабо сходятся к операторнозначным полиномам $P_x(T_x)$, причем все полиномы  $P_x$ различны (динамический индивидуализм).

 Тематика   типичных расширений   обширна, так как  для  каждого инварианта  сохраняющего меру действия заданной группы возникает вопрос о его поднятии (lifting)   типичным расширением. Некоторые общие   задачи пока не решены. Опускается ли типичное свойство автоморфизма на его нетривиальные факторы? Обладает ли типичное расширение   промежуточным расширением? Стабильны ли    лебеговость спектра и свойство кратного перемешивания? 
Аналогичные вопросы возникают при рассмотрении   пространств действий, сохраняющих фиксированную подалгебру и  изучении относительных инвариантов (см. \cite{Sch}).
Интерес   предствляют  типичные расширения  действий на пространствах с сигма-конечной мерой. 
В связи с полными  метриками на пространствах  перемешивающих действий (см. \cite{ST},\cite{Ba}) также возникает  спектр задач о типичных расширениях.

\section{ Примеры типичных свойств расширений}  
Некоторые асимптотические свойства автоморфизма $T$ можно формально определить через  
последовательность значений некоторой функции 
$\varphi(N,j,T)$ при $j\to\infty$. При этом $N$ пробегает натуральный ряд, а 
$\varphi(N,j,T)$ непрерывно зависит  от  $T$.  Напомним примеры таких асимптотических свойств. 

\bf Типичность свойства слабого перемешивания. \rm  
Если для любых множеств $A,B\in\B$ 
выполнено 
$$\mu(A\cap T^j B)\to\ \mu(A)\mu(B),$$
автоморфизм $T$  называется перемешивающим.  Слабое перемешивание означает, что найдется 
перемешивающая последовательность $j_k$, т.е. предыдущая сходимость заменяется на
$\mu(A\cap T^{j_k} B)\to\ \mu(A)\mu(B)$.
Переформулируем эти определения.

Пусть семейство измеримых множеств $\{A_i: i\in\N\}$ плотно в $\B$,  определим функции
$$ \varphi(N,j,T)=\max_{1\leq i,k\leq N} |\mu(A_i\cap T^jA_k)-\mu(A_i)\mu(A_k)|.\eqno (1)$$

Свойство перемешивания эквивалентно тому, что
для любого $N$
$$ \varphi(N,j,T)\to 0, \ \ j\to\infty.$$

Свойство слабого  перемешивания означает, что 
для любых $N$ и $j_0$ найдется $j>j_0$, для которого
$$ \varphi(N,j,T)< \frac 1 N.$$
Это свойство  типично.  Убедимся в этом.
Обозначим через $U_{N,j}$ множество всевозможных $T$, удовлетворяющих последнему неравенству. Оно открыто.
Все слабо перемешивющие автоморфизмы образуют всюду плотное $G_\delta$-множество
$$\bigcap_{N,j_0}\bigcup_{j>j_0} U_{N,j}.$$
Плотность вытекает, например,  из того, что все перемешивающие автоморфизмы плотны в $Aut$ 
(следствие классического факта о всюду плотности  в $Aut$ класса сопряженности апериодического
преобразования).  Халмош доказал, что слабое перемешивание  типично,
а Рохлин установил нетипичность свойства перемешивания. Последнее  вытекает, например,  из 
типичности свойства жесткости, которое несовместимо со свойством перемешивания.

\bf Типичность свойства жесткости. \rm Определим функцию
$$ \psi(N,j,T)= \max \{\mu(A_i)-\mu(A_i\cap T^jA_i)\,:\, 1\leq i \leq N\}.$$
Если для  любых $N$ и $j_0$ найдется $j>j_0$ такое, что 
$$ \psi(N,j,T)< \frac 1 N,$$
преобразование $T$ называется жестким. Обычно это свойство формулируется так: найдется последовательность
$j_k\to\infty$  такая, что $T^{j_k}\to I$ (здесь слабая сходимость операторов совпадает с сильной сходимостью).

Обозначив через $U_{N,j}$ множество всевозможных $T$, удовлетворяющих последнему неравенству, и
положив $\bigcap_{N,j_0}\bigcup_{j>j_0} U_{N,j},$ убеждаемся в том, что все жесткие преобразования
образуют плотное $G_\delta$-множество.
Классическим примером эргодического жесткого преобразования является 
поворот окружности на угол, несоизмеримый с числом $\pi$. 
Мы показали, как устанавливается известный  факт о том, что типичные преобразования являются 
жесткими и слабо перемешивающими. Теперь  приступим к основной теме  статьи.

\bf  Типичные расширения автоморфизма. \rm
 Обозначим через $\J$ семейство всех автоморфизмов пространства $(X\times Y, \mu\otimes\mu')$, 
оставляющих инвариантными  множества вида $A\times Y$ для всех $A\subset X$. Такие автоморфизмы
являются косыми произведениями над  тождественным преобразованием $Id$. 
 Через $Ext(S)$ обозначим все косые произведения $R$  над $S$, их также называют   расширениями.
Напомним, что $R$ определено  формулой
$$R(x,y)= (Sx, R_x y), \ x\in X, \ y\in Y,$$
где $\{R_x\}$ -- измеримое семейство автоморфизмов пространства $(Y,\mu)$. Отметим, 
что мы рассматриваем только случай  $Y=X$, $\mu'=\mu$,  оставляя в стороне конечные расширения,  для которых $|Y|<\infty$. 
Метрика Халмоша на $Aut(\mu\otimes\mu)$ индуцирует полную метрику на замкнутом подпространстве  $Ext(S)$.
Класс  расширений, содержащей $G_\delta$-множество, плотное в $Ext(S)$, называется типичным. Свойство расширения типично, если им обладают представители некоторого типичного класса. 
Для доказательства типичности ряда свойств применяется следующее хорошо известное утверждение, которое несложно установить при помощи  классической леммы Рохлина-Халмоша:
\it  для всякого  косого произведения $R$ над эргодическим автоморфизмом $S$  класс
 $\{\Phi^{-1}R\Phi :  \Phi\in \J\}$ плотен в  $Ext(S)$. \rm
В частности, это утверждение леммы выполняется для прямых произведений  $R=S\times T$, где $T$ -- некоторый автоморфизм.

\bf Поднятие свойства   частичной жесткости. \rm 
Фиксируя параметр $a\in (0,1]$ и плотное в $\B$ семейство $\{A_i\}$,  определим функцию
$$ \psi_a(N,j,T)= \max \left\{a\mu(A_i)-\mu(A_i\cap T^jA_i)\,:\, 1\leq i \leq N\right\}.$$
 
Если для  любых $N$ и $j_0$ найдется $j>j_0$ такое, что 
$$ \psi_a(N,j,T)< \frac 1 N,$$
то преобразование $T$ называется $a$-жестким. Иначе это свойство формулируется так: найдется последовательность
$j_k\to\infty$  такая, что $T^{j_k}\to_w aI+(1-a)P$, где $P$ -- некоторый марковский оператор
(марковость оператора означает, что $P$ и $P^\ast$ сохраняют неотрицательность функций, а константы переводят в себя).  

\vspace{3mm} \bf
Теорема 2.1.    \it  Типичные расширения эргодического $a$-жесткого преобразования сохраняют свойства  $a$-жесткости.
 \rm

\vspace{3mm} 
Доказательство. 
 Расширение  $R$ преобразования $S$ обладает $a$-жесткостью, 
если для  любых $N$ и $j_0$ найдется $j>j_0$ такое, что
$\psi_a(N,j,R)< \frac 1 N$. Тогда $G_\delta$-множество
$$ W_a=\bigcap_N\bigcap_{j_0} \bigcup_{j>j_0}\left\{R:\psi_a(N,j,R)< \frac 1 N\right\}$$
состоит в точности из $a$-жестких преобразований. Так как произведение 
 $S\times Id$ и все косые произведения вида  $J^{-1}(S\times Id)J$, $ J\in\J$, наследуют  
свойство  $a$-жесткости, $W_a$ всюду плотно в $Ext(S)$. 
Таким образом, множество $W_a$ всех $a$-жестких преобразований  типично.

\bf Поднятие  свойства слабого перемешивания. \rm Для  слабо перемешивающего  автоморфизма $S$ 
рассматриваем расширения $R$ и функцию $ \varphi(N,j,R)$, определенную формулой (1) с оговоркой, 
что теперь вместо меры $\mu$ фигурирует мера $\mu\otimes\mu$. Если в предыдущем доказательстве 
заменить $ \psi_a(N,j,R)$ на 
$\varphi(N,j,R)$, а вместо $S\times Id$ рассмотреть слабо перемешивающее расширение $S\times S$,  
аналогично получим, что слабо перемешивающие расширения образуют плотное  $G_\delta$-множество.

Для косого произведения $R=(S,R_x)$ над эргодическим преобразования $S$ определено понятие 
относительного слабого перемешивания, которое означает эргодичность косого произведения
$R\times_S R:=(S,R_x\times R_x)$ относительно меры $\mu\otimes\mu\otimes\mu$. Пишем в этом случае
$R\in RWM(S)$. В общем случае (если  не требовать эргодичность  $S$) относительное слабое перемешивание
можно определить так: для любых $A,B\in\B$ выполнено
 $$ \int_X \frac 1 j \sum_{n=1}^j(\mu(C(x,n, R)A\cap B)-\mu(A)\mu(B))^2 d\mu(x)\ \to 0,\ j\to\infty,$$
где $$C(x,n,R)=R_{S^{n-1}x}\dots R_{Sx}R_x.$$
Это свойство можно переформулировать: для больших значений $N$ верно, что в среднем по $x$ и в среднем по $n$, $1\leq n\leq N$, операторы $C(x,n, R)$ близки к $\Theta$ в слабой операторной топологии.
Напомним, что $\Theta$ означает ортопроекцию на пространство констант.
Приведем вместе с доказательством результат Глазнера-Вейса \cite{GW}.

\vspace{3mm}
\bf Теорема 2.2.     \it  Если $R$ -- типичное расширение эргодического автоморфизма $S$, то  $R \in RWM(S)$.  \rm

\vspace{3mm}
 Доказательство. Пусть семейство измеримых множеств $\{A_i: i\in\N\}$ плотно в $\B$, положим
$$ \varphi(N,j,R)=\max_{1\leq i,k\leq N} \int_X \frac 1 j \sum_{n=1}^j\left(\mu(C(x,n, R)A_i\cap A_k)-\mu(A_i)\mu(A_k)\right)^2 d\mu(x).$$

Условие $R\in RWM(S)$ означает, что   для любого натурального $N$ найдется $j$ такое, что $ \varphi(N,j,R)<\frac 1 N$.
Так $\varphi$ непрерывно зависит от $R$, класс $RWM(S)$ является $G_\delta$-множеством. Плотность класса $RWM(S)$ в пространстве  $Ext(S)$  вытекает из плотности
его подкласса $\{J^{-1}(S\times T)J:\ J\in\J\}$, где $T$ --  некоторый слабо перемешивающий автоморфизм. Теорема доказана.

\section{Слабые замыкания и динамические свойства типичных семейств }
Функцию $P$ от унитарного оператора $T$ назовем допустимой, если 
$P(T)=c\Theta +\sum_{i=0}^\infty c_i T^i$,
где $c,c_i\geq 0$, $c + \sum_{i=0}^\infty c_i=1$.

\bf Допустимые слабые пределы. \rm  Если имеет место слабая сходимость $S^{n_i}\to P(S)$, где $P$ -- допустимая функция, 
будет ли оператор  $P(R)$  слабым пределом  степеней  типичного расширения $R$?
Ответ положительный для всевозможных  факторов $S$ типичных преобразований $T$.  

\vspace{1mm} \bf
Теорема 3.1.    \it Если степени  некоторого расширения $\tilde R$ преобразования  $S$ обладают слабым пределом  $P(\tilde  R)$, то множество всех расширений $R$ с таким  свойством является типичным.
\rm

Доказательство. Рассмотрим функцию 
$$\varphi(N,j,R)= \max_{1\leq m,n\leq N} \left((R^j-P(R))f_m,\,f_n\right),$$
где $\{f_n\}$ -- плотное в единичном шаре семейство функций.
Положим 
$$ W=\bigcap_N\bigcap_{j_0} \bigcup_{j>j_0}\left\{R:\varphi(N,j,R)< \frac 1 N\right\}.$$

 Множество $W$ типично, так как содержит всюду плотный в $Ext(S)$ класс расширений, 
когомологичных расширению $\tilde R$, т. е  сопряженных с $\tilde R$ посредством 
некоторого $\Phi\in\J$. В то же время $W$ содержит все $R$ такие, что 
слабое замыкание их  степеней содержит оператор $P$.

Замечание.  Теорему можно усилить, если  в формулировке  заменить одну функцию  $P$ на произвольное семейство допустимых функций. Это непосредственно  вытекает из того, что множество всех операторов вида $P(R)$,  где $P$ пробегает заданное  семейство  допустимых функций, является сепарабельным в слабой операторной топологии.

\bf Типичные семейства. \rm Косому произведению $(Id, T_x)$ над тождественным оператором 
отвечает  измеримое семейство автоморфизмов $\{T_x: x\in X\}$.  Типичному множеству в пространстве $Ext(Id)$ соответствует типичное множество таких семейств. Из теоремы 3.1  вытекает, что для типичного семейства почти все входящие в него автоморфизмы обладают слабым перемешиванием. Ниже мы установим другие типичные свойства измеримых семейств, выражаемые в терминах  слабых пределов. 
Будет показано, частности, что для типичного семейства входящие в него автоморфизмы при итерациях  на некоторых последовательностях  сходятся к тождественному оператору, хотя  на некоторых 
других последовательностях одна часть из них сходится к тождественному оператору $I$, а другая часть сходится, например, к оператору $\Theta$. 
Приведем пример. Для удобства считаем, что $X=[0,1]$.  
Для типичного семейства $\{T_x: x\in [0,1]\}$ верно, что найдется такая последовательность $n_i$, что для почти всех $x$ выполняется
$$T_x^{n_i}\to_w\ x \Theta+(1-x)I.$$ 
Теперь рассмотрим  более общую ситуацию.  Пусть $c(x),c_i(x)$ -- неотрицательные измеримые функции, 
удовлетворяющие условию $c(x) + \sum_{i=0}^\infty c_i(x)=1$.  Для типичного семейства $\{T_x: x\in [0,1]\}$ найдется  последовательность
$n_i$ для которой 
$$T_x^{n_i}\to_w\ c(x)\Theta +\sum_{i=0}^\infty c_i(x) T^i_x.$$
Сказанное является следствием   следующего утверждения.

\vspace{2mm} \bf
Теорема 3.2.    \it  Пусть $X=\uu_m A_m$,  все множества   $A_m$ имеют положительную меру, а также  задан   
 набор допустимых функций $\{P_m: m\in\N\}$ и последовательности $n_j\to\infty$. Тогда 
  для типичного расширения  тождественного оператора
$(Id,T_x)$ найдется последовательность $n_{j(k)}\to\infty$ такая, что для каждого $m$ для почти всех
$x\in A_m$ имеют место  слабые сходимости $T_x^{{j(k)}}\to P_m(T_x)$.\rm 

\vspace{2mm} 
Доказательство.  Из теоремы 2.1 \cite{R21} вытекает существование счетного семейства автоморфизмов $U_m$ таких, что для некоторой последовательности  $n_j\to\infty$ выполнено $U^{n_j}\to P_m(U)$ при $j\to\infty$.
    Пусть $R=(Id, R_x)$, где  $R_x=U_m$   при $x\in A_m$.  Класс расширений $\{J^{-1}RJ:J\in\J\}$ 
 плотен в $Ext(Id)$ и лежит в множестве расширений $(Id, \tilde R_x)$, для которых для всех $m$ выполнено 
$\tilde R_x^{n_j}\to  P_m(\tilde R_x)$ при $x\in A_m$. 
Стандартные рассуждения (похожие на рассуждения в доказательстве теоремы 3.1)  показывают, что множество  
всех расширений $W$, для которых такая сходимость выполняется для какой-нибудь  последовательности $j(k)\to\infty$, является $G_\delta$-множеством. Теорема доказана.

\section{Типичные  расширения сохраняют сингулярность спектра}
  
Если мера $\sigma$  на единичной окружности $\T$ в комплексной плоскости
  сингулярна относительно меры Лебега $m$, то это эквивалентно свойству 

 $(\ast)$: \it  для любого $N>0$ существует натуральное $P$ такое, что для разбиения окружности на   дуги $I_{k,P}=\left[e^{2\pi i \frac k P}, e^{2\pi i\frac {k+1} P}\right]$, 
$k=0, 1,\dots, P-1$, выполнено неравенство
  $$\left|\left\{k: \sigma(I_{k,P})<1/ NP\right\}\right|>\left(1-\frac 1 N\right)P.$$ \rm

Поясним, почему  $(\ast)$ влечет за собой сингулярность меры. Если мера $\sigma$ имеет абсолютно непрерывную компоненту $\nu$, то для некоторых $a,b>0$  найдется множество $A$ меры $2a$, на котором производная Радона-Никодима меры $\nu$ больше $2b$.
Для больших $P$   наилучшим образом приблизим множество $A$  объединениями дуг $I_{k,P}$.  Большинство этих дуг состоит в основном их точек множества $A$, суммарная   мера этого большинства  больше  $a$.     Получили, что для достаточно большого $P$ для более, чем $aP$ дуг верно, что они  имеют $\sigma$-меру больше $b/P$.
В $(\ast)$  этот случай запрещен при $N$, когда $1/N$    меньше  $a$ и $b$.
Таким образом, из $(\ast)$ получили сингулярность меры.

Покажем, что из сингулярности меры  $\sigma$  вытекает $(\ast)$. Пусть $\sigma(\T)=1$. Из общих фактов вытекает, что для всяких  $\eps>0$ и $N$ найдется $P'$ и 
конечное объединение $U$ некоторых  дуг $I_{k,P'}$ таких, что  $m(U)<1/3N$, причем  
 $\sigma(U)>1- \eps/2$.    Для всех достаточно больших $P>P'$  в дополнении к $U$ найдется множество $V$, состоящее из некоторых дуг $I_{k,P}$ ($k$ меняется, $P$ фиксировано) таких, что $m(V)>(1-1/2N)$ и  $\sigma(V)<\eps$.  Пусть $\eps < 1/2N^2$, тогда для достаточно больших $P$ найдется не более   $P/2N$ дуг $I_{k,P}$ из $V$ с $\sigma$-мерой больше  $1/NP$ (иначе 
$\eps > 1/2N^2$). 
Дополнение к $V$ состоит из не более, чем  $P/2N$ дуг. Следовательно, число дуг, для которых их $\sigma$-мера меньше $1/NP$,  превосходит $P-P/N$ для всех достаточно больших $P$.
Это означает выполнение свойства $(\ast)$.

Переформулируем    свойство $(\ast)$ в терминах  непрерывных функций на окружности. Зададим  непрерывные функции $\Delta_{k,P}$  на окружности $\T$ следующим образом:   
на дуге $I_{k-1,P}$, $k\in\Z_P$, функция $\Delta_{k,P}$ растет линейно от 0 до 1, на  $I_{k,P}$ тождественно равна 1,  на дуге   $I_{k+1,P}$  линейно убывет до  0 и принимает нулевое значение на  остальных дугах.
Положим
  $$D(\sigma,N,P)=\left\{k: \int_\T\Delta_{k,P} d\sigma < \frac 1 {NP}\right\}.$$
Мера $\sigma$  на окружности $\T$  сингулярна только в том случае, когда
 для любого $N$ найдется $P$ такое, что
$$|D(\sigma,N,P)|>\left(1-\frac 1 N\right)P.$$

\vspace{3mm}
\bf Теорема 4.1.    \it  Типичные расширения сохраняют сингулярность спектра.
\rm

\vspace{3mm}
Доказательство.  
Пусть $\sigma_{f,R}$ обозначает спектральную меру оператора $R$ с цикличесим вектором $f$, $\|f\|=1$, т.е.
$$\widehat{\sigma_{f,R}}(s)=\int_\T z^s d\sigma_{f,R}=(R^sf,f).$$
Непрерывная функция $\Delta_{k,P}$ равномерно близка к некоторой  сумме Фейера $S_n$,
значение $\int_\T\Delta_{k,P} d\sigma_{f,R}$ близко к  интегралу $\int_\T S_n d\sigma_{f,R}$,
который  непрерывно зависит лишь от конечного набора коэффициентов Фурье  меры $\sigma_{f,R}$.
Так как коэффициенты $\widehat{\sigma_{f,R}}$ непрерывно зависят от $R$, получаем, что 
 множество  
 $$\left\{R: \int_\T\Delta_{k,P} d\sigma_{f,R} < \frac 1 {NP}\right\}$$
открыто. Тогда 
$\{R\,: \,|D(\sigma_{f,R},N,P)|= d\}$ открыто, значит, открыто множество 
$$U(f,N,P)=\left\{R\,:\, |D(\sigma_{f,R},N,P)|>   {\left(1-\frac 1 N\right)P}  \right\}.$$
Из сказанного вытекает, что множество автоморфизмов $R\in Aut(\mu\otimes\mu)$ с сингулярной мерой 
$\sigma_{f,R}$ является $G_\delta$-множеством
$$Sing(f)=\bigcap_N\bigcup_P U(f,N,P).$$ 
В пространстве $L_2$ выберем ортонормированный базис $\{f_i\}$ и положим 
$$Sing=\bigcap_i Sing(f_i).$$
Мы доказали, что автоморфизмы $Aut(\mu\otimes\mu)$ с сингулярным спектром образуют $G_\delta$-множество. Косые произведения
над  $S$ образуют замкнутое множество  $Ext(S)\subset Aut(\mu\otimes\mu)$, а расширения с сингулярным спектром   образуют $G_\delta$-множество $Sing_S$ в индуцированной топологии на $Ext(S)$.
Осталось  заметить, что  $Sing_S$ содержит  всюду плотное множество
 $\{J^{-1}(S\times Id)J :  J\in \J\},$  где  $\J$ обозначает класс косых произведений над $Id$. Теорема доказана.

 В теореме  сингулярность спектра не запрещала наличие дискретной компоненты. Если $S$ имеет непрерывный  сингулярный спектр, то его типичное расширение наследует  это  свойство.   

\bf Замечания о спектральных кратностях. \rm Сохранение  некоторых асимптотических свойств (инвариантов)  для типичных расширений было доказано методом, в котором принципиальную  роль играет пример косого произведения, обладающего заданным  асимптотическим свойством.  В случае  частичной жесткости и сингулярного спектра подходящим примером было произведение $S\times Id$. Для слабого перемешивания расcматривалось произведение $S\times S$. Следующая теорема дает примеры расширений, сохраняющих наборы спектральных кратностей базового автоморфизма.

\vspace{2mm}
\bf Теорема 4.2. \it  Если если для автоморфизма $S$ выполнено $S^{n_j}\to I$, $n_j\to\infty$,   то для типичного автоморфизма $T$ произведение $R=S\times T$  имеет набор спектральных кратностей такой же, как у $S$.
 \rm

\vspace{2mm}  
Доказательство. Пусть  $S$ имеет простой спектр. Среди типичных преобразований найдется  $T$ с простым спектром и свойством  $T^{n_{j_k}}\to T$  для некоторой последовательности 
$j_k\to\infty$ (вытекает из   \cite{R21}, теорема 2.1).
Тогда произведение $S\otimes T$ имеет простой спектр, так как тензорное произведение циклических векторов $f$ и $g$ операторов $S$ и $T$, является циклическим вектором для оператора $S\otimes T$.
Действительно, для оператора $S\otimes T$ циклическому пространству $C$, содержащему вектор $f\otimes g$ будут принадлежать все векторы вида $f\otimes T^ng$. Следовательно, ему будут принадлежать
все векторы вида $S^mf\otimes T^{m+n}g$, значит, $C$ совпадает со всем пространством $L_2\otimes L_2$.
 
Если $S$ является прямой суммой набора операторов $S_i$ с простым непрерывным спектром, $S^{n_j}\to I$, то  $S\otimes T$ являтся прямой суммой операторов $S_i\otimes T$, которые по тем же причинам, что    выше, имеют простой  спектр. Если $S_i$ и $S_{i'}$ имеют взаимно сигулярные  спектральные типы, то операторы
$S_i\otimes I$ и $S_{i'}\otimes I$ также имеют взаимно сигулярные  спектральные типы, следовательно, они не допускают 
ненулевое сплетение. Но тогда нет ненулевого сплетения между  $S_i\otimes T$ и  $S_{i'}\otimes T$.
Покажем это. Пусть для некторого оператора $U$ выполнено
$$(S_i\otimes T)U=U(S_{i'}\otimes T),$$
тогда
$$(S_i\otimes T)^{-n_{j_k+1}}U=U(S_{i'}\otimes T)^{-n_{j_k}+1},$$
$$(S_i\times I)U=U(S_{i'}\times I), \ \ U=0.$$ 
Таким образом,  набор спектральных кратностей у $S_i\otimes T$ такой же, как у оператора $S$. Теорема доказана.

Без доказательства сформулируем следующее утверждение: \it если $S^{n_j}\to aI+ (1-a)\Theta$, $a\in (0,1]$, то типичное расширение автоморфизма $S$ сохраняет кратности спектра.\rm


\section{Асимметрия, кратные слабые пределы}
 Если преобразования $R$ и  $R^{-1}$ не сопряжены в $Aut$,  то 
такое $R$ называем  асимметричным.
Типичное расширение симметричного преобразования с положительной энтропией,
сохраняет свойство симметричности (следствие результата \cite{AGTW}).
Мы докажем, что для жестких преобразований типичное расширение асимметрично. 

Следующее вспомогательное для наших целей утверждение  является несложной модификацией  основного результата   работы  \cite{R03} (последовательность $m_j\to\infty$ теперь играет роль  подпоследовательности наперед заданной  последовательности).

\bf  Теорема 5.1.  \it  Существует    преобразование $T$ со следующим свойством:  для  любой последовательности, стремящейся к бесконечности, найдется ее подпоследовательность  
$m_j\to\infty$ такая, что   для любого измеримого множества $A$ имеют место сходимости
$$4\mu(A\cap T^{m_j}A\cap T^{3m_j}A)\to \  \mu(A)+ \mu(A)^2 + 2\mu(A)^3,$$
$$\mu(A\cap T^{-m_j}A\cap T^{-3m_j}A)\to \ \mu(A)^2.$$

\bf  Теорема 5.2. \it Типичное расширение жесткого  эргодического преобразования является асимметричным.
\rm

\vspace{3mm}
Доказательство.  Рассмотрим $R=S\times T$, $\bar \mu=\mu\times \mu$, где  $T$ удовлетворяет условиям теоремы 5.1, а $S$ -- жесткое преобразование, причем $S^{m_j}\to I$.
Тогда для любого $\bar A$, $\bar \mu(\bar A)>0$ выполняется
$$\lim_j \bar\mu( \bar A\cap R^{m_j}\bar A \cap R^{3m_j}\bar A )
>  \frac 1 3  \bar\mu(\bar A).$$
Для   множества $\bar A_0=X \times A_0$, $\bar \mu(\bar A_0)=\frac 1 4$, имеем
$$\lim_j \bar\mu( \bar A_0\cap (R^{-m_j}\bar A_0 \cap R^{-3m_j}\bar A_0 )=\frac 1 {16}.$$

Рассмотрим семейство $W$ всех $V\in Ext(S)$, удовлетворяющих условию: для любых $i$, $j_0$ найдется $j>j_0$  такое, что  
$$\bar\mu(\bar A\cap V^{m_j}\bar A \cap V^{3m_j}\bar A )>  \frac 1 3 \bar\mu(\bar A) $$
и
$$ \bar\mu( \bar A_0\cap V^{-m_j}\bar A_0 \cap V^{-3m_j}\bar A_0 )< \frac 1 {15}.$$
Семейство  $W$ является $G_\delta$-множеством. Оно содержит все косые произведения, когомологичные расширению $R$, т.е. класс $\{J^{-1}RJ: J\in\J\}$, поэтому является всюду плотным в $Ext(S)$. 

Остается заметить, что $W$ состоит из несимметричных автоморфизмов.
Действительно, пусть $\Phi^{-1}V^{-1} \Phi = V\in W.$
Перепишем последнее неравенство в виде
$$ \frac 1 {15} > \bar\mu( \bar \Phi^{-1} \Phi A_0\cap V^{-m_j}\Phi^{-1}\Phi \bar A_0 
\cap V^{-3m_j}\Phi^{-1}\Phi\bar A_0 ).$$
Для больших $j$ имеем
$$ \frac 1 {15} >  \bar\mu( \Phi \bar A_0\cap  \Phi V^{-m_j}\Phi^{-1}\Phi \bar A_0 
\cap  \Phi V^{-3m_j}\Phi^{-1}\Phi\bar A_0 )= $$
$$
=\bar\mu( \Phi \bar A_0\cap  V^{m_j}\Phi \bar A_0 
\cap  V^{3m_j}\Phi\bar A_0 ) > \frac 1 {12}.$$
Полученое противоречие показывает, что $V$, $V^{-1}$ неизоморфны.
Теорема доказана.

Таким образом,  возможность поднять  свойство симметричности при типичном расширении, как мы увидели,
 зависит от перемешивающих свойств базового автоморфизма.
 Представляется весьма правдоподобной гипотеза о том, что симметричность не сохраняется при типичных расширениях всех частично жестких и  некоторых перемешивающих преобразований.
Если симметричное эргодическое преобразование имеет положительную энтропию, то из результата 
\cite{AGTW} непосредственно вытекает, что его типичное расширение сохраняет симметричность. 
В связи с этим возникает вопрос:  существуют ли симметричное преобразование  с нулевой энтропией, для которого  
его  типичные расширения также симметричны?  

\bf Замечание.  \rm Имеются преобразования $T$  с необычным  нетипичным свойством: все  декартовы 
конечные степени $T\times\dots\times T$ асимметричны, но бесконечная декартова степень $T\times T\times \dots$ симметрична. Указание: $T=S^{-1}\times S\times S$.


\section{ Типичные расширения имеют бесконечную  $P$-энтропию }
 Типичное преобразование,  неформально говоря, обязано хорошо перемешивать на чрезвычайно  длинных временных интервалах. Например,  перемешивание имеет место  на интервалах вида 
$$\Large \left(n_j, \ n_j^{n_j}!\right)$$ 
для некоторой последовательности $n_j\to\infty$. Следующее   утверждение (i) является частным случаем теоремы 2 из \cite{R}. Пункт (ii) доказывается аналогично пункту (i).
 
\vspace{3mm}
 \bf Теорема 6.1. \it  Для любой последовательности конечных множеств $P_j\subset \N$, удаляющихся от 0, 

(i) для типичного преобразования $T$  найдется подпоследовательность $j(k)\to\infty$ такая, что $P_{j(k)}$ является перемешивающей последовательностью;

(ii)   для типичного расширения  $R$ перемешивающего преобразования $S$ найдется подпоследовательность $j(k)\to\infty$ такая, что $P_{j(k)}$ является перемешивающей последовательностью.
\rm 
\vspace{3mm} 

Утверждение (i) получило развитие в терминах энтропийных инвариантов в работе \cite{R21}. Метод этой работы 
ниже применяется для  типичных расширений. В частности, будет получено новое доказательство  результа \cite{AGTW} о том, что типичное расширение системы $S$ с нулевой энтропией неизоморфно $S$. Мы подберем энтропийный инвариант, который равен 0 для такого $S$, но принимает бесконечное значение для типичного расширения преобразования $S$.

\vspace{2mm} 
\bf $P$-энтропия Кушниренко. \rm 
Определение   $P$-энтропии автоморфизма $T$ из \cite{R21} является удобной для нас 
  модификацией   определения последовательной энтропии Кушниренко \cite{Ku}.  Для последовательности $P$ конечных множеств $P_j\subset \N$  и автоморфизма $T$ вероятностного пространства $(\bar X,\bar\mu)$ определим энтропию $h_P(T)$ следующим образом. Положим 
 $$h_j(T,\xi)=\frac 1 {|P_j|}  H\left(\bigvee_{p\in P_j}T^p\xi\right),$$
где 
$\xi=\{C_1,C_2,\dots, C_n\}$ -- измеримое разбиение множества $X$. Напомним, что энтропия разбиения
определяется формулой 
$$ H(\xi)=-\sum_{i=1}^n \bar\mu( C_i)\ln \bar\mu( C_i).$$
Теперь положим
$$h_{P}(T,\xi)={\limsup_j} \ h_j(T,\xi),$$

$$h_{P}(T)=\sup_\xi h_{P}(T,\xi).$$

 Типичные преобразования обладают плохими перемешивающими свойствами, к чему давно привыкли специалисты. Однако, на некоторых последовательностях они могут перемешивать лучше, чем, например,   автоморфизм, входящий в орициклический поток $O_t$. 
Пусть $$P_j=\{2^i:  n(j)\leq i < n(j+1)\}, \ \  {n(j+1)/{n(j)}}\to \infty,$$
Для автоморфизма  $T=O_1$, входящего в орициклический поток $O_t$, 
из результатов работы \cite{Ku}  вытекает, что  $0<h_P(T)<\infty$.
Как известно,  $T$ обладает кратным перемешиванием и лебеговским спектром.
Типичный автоморфизм $S$  вероятностного пространства имеет сингулярный спектр,
более того $S$ является жестким (т.е. $S^{n_i}\to I$ для некоторой  последовательности     
${n_i}\to \infty$). Однако,   для типичного автоморфизма выполнено  $h_P(S)=\infty$,
как показано в \cite{R21}.  Получается, что на некоторой подпоследовательности 
$P_{j(k)}$ типичный автоморфизм $S$ перемешивает гораздо  лучше, чем $T$.

Ниже для удобства мы ограничимся рассмотрением частного  случая, когда
 $P_j$ -- последовательность  расширяющихся арифметических прогрессий.

\vspace{3mm}\bf 
 Лемма (\cite{R21}). \it Если  $h(S)=0$, то для некоторой последовательности 
$$P_j=\{j,2j,\dots, L(j)j\}, \ L_j\to\infty,$$ выполнено  $h_P=0$. 
Если  $h(T)>0$, то   $h_P(S\times T)=\infty$. \rm

\vspace{3mm}
Первое утверждение   вытекает из того, что для любого $j>0$ энтропия $h(S^j)$ равна 0.
Второе утверждение  леммы следует   из того, что  $h(T^j)=jh(T)$.

\vspace{2mm} \bf
Теорема 6.2.    \it  Множество  $\{R\in Ext(S): h_P(R)=\infty\}$ типично. \rm

\vspace{2mm}
Доказательство. Фиксируем семейство автоморфизмов    $\{J_q\}$, $q\in \N$, плотное в  $\J$.
Положим  $R_q=J_q^{-1}RJ_q$, где $R=S\times T$, а $T$ -- бернуллиевский автоморфизм с образующим разбиением
$\{C_1, C_2,\dots, C_k\}$.
  Семейство   $\{R_q\}$ плотно в   $Ext(S)$.
Рассмотрим  разбиение   $\xi$, $H(\xi)>0$,  множества  $X\times Y$ вида
$$\xi=\{X\times C_1, X\times C_2,\dots, X\times C_k\}.$$

Покажем, что для больших значений $j$ выполнено 
$$ h_{j}(R_q,\xi)=h_{j}(R,J_q\xi) =\frac 1 {L_j} H\left(\bigvee_{n=1}^{L(j)} R^{nj}J_q\xi\right) >H(\xi)/2.\eqno (2) $$
Косое произведение $J_q$  (номер $q$ фиксирован) имеет вид $J_q(x,y)=(x, Q_x y).$
Измеримое семейство автоморфизмов $\{Q_x\}$, $x\in X$, как функция от $x$
 приближется по мере $\mu$ на $X$ конечнозначной измеримой функцией, принимающей значения
из некоторого конечного множества автоморфизмов $\{\tilde Q_x: x\in X\}$. Здесь подразумевается приближение автоморфизмов $Q_x$ автоморфизмами $\tilde Q_x$  относительно метрики Халмоша. 
Разбиения вида $\{\tilde Q_x C_1, \dots, \tilde Q_x C_k\}$  образуют конечное множество, обозначим эти разбиения через
$\Delta_d$, $1\leq d\leq D$. Каждое $\Delta_d$  приближаются разбиениями, измеримыми относительно 
$$\eta_M=\bigvee_{i=-M}^M \{T^iC_1, \dots, T^iC_k\}$$
для некоторого достаточно большого натупального числа $M$.
Замечаем, что  при  $j>2M$ разбиения 
$T^{nj}\eta_M$, $n=1, 2, \dots,$ независимы. Это означает почти независимость 
разбиений $\{T^{nj} \Delta_{d_n}\}$ при любом выборе последовательности $d_n$, $1\leq d_n\leq D$, и как следствие  почти независимость разбиений $R^{nj}J_q\xi$,  $n\in \N$, которая обеспечивает выполнение неравенства  (2).

Таким образом,  для любых   $q$, $N$ найдется   $j=j(q,N)$ и окрестность  
$U(q,N)$ автоморфизма  $R_q$ такие, что  
$j>N$ и для всех  $V\in U(q, N)$ выполнено неравенство
$$ h_{j}(V,\xi) > H(\xi)/2.$$
Множество
$$W=\bigcap_N\bigcup_q U(q,N)$$ является плотным $G_\delta$-множеством.  Действительно, если  $ V\in W$,
то для каждого $N$ найдется   $q(N)$ такое, что
$$ h_{j(q(N),N)}(V,\xi)> H(\xi)/2.$$ 
Так как   $j(q(N),N)>N$, очевидным образом получим $$h_P(V)\geq H(\xi)/2.$$  
В силу того, что бернуллиевские разбиения $\xi_i$ можно выбрать с условием  $H(\xi_i)\to\infty$, 
а пересечение сооветствующих типичных множеств $W_i$
типично, получаем  $ h_P(V)=\infty$ для всех $V\in\bigcap_i W_i$. 
 Теорема доказана.

\bf Тонкие семейства. \rm
По контрасту со свойством доминантности, обнаруженным в \cite{AGTW} для систем с положительной энтропией, будем называть семейство $F$ автоморфизмов тонким  (exquisite), если для каждого $S\in F$   типичное множество его расширений  не содержит представителя, изоморфного какому-нибудь  $S'\in F$. 
 Из  теоремы 6.2  вытекает,  что  \it  множество  $\{S: h_P(S)<\infty\}$ является тонким. \rm

\section{Рекуррентность типичных коциклов,  отсутствие  независимого фактора и стабильность перемешивания } 

 Говорим, что свойство  действия  стабильно, если типичное расширение системы с этим свойством таже им обладает. Непрерывность спектра,   сингулярность спектра, частичная  жесткость,  детерминированность системы, К-свойство, как известно,  стабильны. К этому списку сейчас  мы добавим свойства  мягкого и строгого  перемешивания.  Открытыми, напомним,  остались вопросы о стабильности свойств автоморфизма иметь лебеговский  спектр и обладать кратным перемешиванием.

\vspace{2mm}
 \bf Теорема 7.1. \it  Свойство перемешивания стабильно.
\rm 
\vspace{2mm} 

Доказательство этой теоремы использует  вспомогательные утверждения, которые представляют самостоятельный интерес.

\bf Неперемешивающий фактор,  независимый со  всяким перемешивающим 
фактором. \rm
 Нам понадобится теорема Ф.~Парро, полученная им   в  2002 году. 

\vspace{1mm}  
\bf Теорема (\cite{P}). \it  Неперемешивающий эргодический автоморфизм обладает нетривиальным фактором, который дизъюнктен со всеми перемешивающими автоморфизмами. 

\vspace{1mm}  
\rm Дизъюнктность влечет за собой независимость фактора Парро от  всякого  перемешивающего фактора автоморфизма $R$, чем мы воспользуемся ниже.

Приступим к доказательству теоремы 7.1.  Пусть расширение $R$ перемешивающего автоморфизма 
$S$ не является перемешивающим, тогда оно обладает фактором Парро, который независим от базового $S$-фактора.   Напомним, что фактором называется ограничение действия на инвариантную сигма-подалгебру.
Стабильность свойства  перемешивания вытекает из сказанного и  следующего утверждения. 

\vspace{2mm}
 \bf Теорема 7.2. \it  Типичное расширение $R$ слабо перемешивающего действия $S$ не обладает нетривиальным фактором,  независимым от базового $S$-фактора.
\rm 
\vspace{2mm}

Для последовательности конечных множеств $P_j\subset\N$, расширения $R$ и множества $A\subset X$, $\mu(A)>0$,
определим 
$$\varphi_A(N,j,R)=\prod_{p\in P_j}\mu\left(x\in A\cap S^pA: \rho(C(x,p,R), Id)<\frac 1 N\right),$$
где $\rho$ -- метрика Халмоша в $Aut(\mu)$, а  $C(x,n,R)$ -- обозначение для коцикла $R_{S^{n-1}x}\dots R_{Sx}R_x.$  
 Будем писать $R\in RC(P_j,A)$,  
если для каждого $N$ и $j_0$ найдется такой  $j>j_0$, что  $\varphi_A(N,j,R)>0$. 

Следующее утверждение можно назвать теоремой о рекуррентности коцикла $C(x,n,R)$, отвечающего типичному расширению $R\in Ext(S)$.

\vspace{2mm}
 \bf Теорема 7.3. \it  Для каждого множества $A$ положительной меры  
типичное расширение $R$ фиксированного слабо перемешивающего преобразования $S$  принадлежит классу $RC(P_j,A)$ для некоторой последовательности 
 $P_j\subset \{j,j+1,\dots,2j\}$, причем $|P_j|/j\to 1$ при $j\to\infty$. \rm 

\bf Следствие. \it  Для слабо перемешивающего преобразования  $S$, его типичного  расширения $R=(S,R_x)$ и множества $A$ положительной меры  для почти всех $x\in A$ найдется последовательность $p_i\to\infty$ (она зависит от $x$)
такая, что $C(x,p_i,R)\to Id$, причем  $S^{p_i}(x)\in A$. \rm

\vspace{2mm}
Доказательство теоремы 7.3. Множество $P_j$ выбираем так, чтобы   $S^{p(j)}\to\Theta$ для всякой последовательности
 $p(j)\to\infty$ при условии  $p(j)\in P_j$.
Так как слабое перемешивание эквивалентно перемешиванию на множестве плотности 1, дополнительно обеспечиваем условие 
 $|P_j|/j\to 1$ при $j\to\infty$. Класс  $RC(P_j,A)$ является $G_\delta$-множеством, что вытекает из его определения и
непрерывности зависимости $\varphi_A(N,j,R)$ от $R$. 
Докажем  его плотность.  Для этого рассмотрим класс расширений, когомологичный тривиальному расширению $S\times Id$. Заметим, что 
$$C(x,p,J^{-1}R_0J)=J^{-1}_{S^p(x)}J_x.$$
Для любого $\eps>0$  найдется множество $A'\subset A$ такое, что $\mu(A')>0$ и для некоторого $x_0\in A'$ и всех $x\in A'$ выполнено
 $$\rho(J_x, J_{x_0})<\eps.   $$ 
Тогда при заданном $N$ выбор достаточно малого числа $\eps$ обеспечивает для всех $x\in A'$ при условии $S^p(x)\in A'$  выполнение неравенства 
 $$ \rho(J^{-1}_{S^p(x)}J_x, Id)<\frac 1 N.\eqno(3)$$
Но в силу перемешивающих свойств степеней $S^p$ для всех достаточно больших $p\in P_j$ множество таких $x\in A'$, что  $S^p(x)\in A'$,  имеет положительную меру.
Теорема доказана.

Доказательство теоремы 7.2. Пусть множество $E \subset X\times Y$ независимо от $S$-фактора,
это означает, что функция $h(x)=\mu(y: (x,y)\in E)$  п.в. равна $e$ -- мере множества $E$.  Если $E$ принадлежит
$R$-инвариантной алгебре, независимой от $S$-фактора, то  $h_p(x)=\mu(y: (x,y)\in E\cap R^pE)$
также является п.в. константой $h_p$, причем    в силу перемешивания
$R^p$   при  $p\in P_j$, $p\to\infty$, выполнено $h_p\to e^2$.   Из  (3) мы видим, что для множества $x$-ов положительной меры $h_p(x)$ близко к $e$.   Но $h_p(x)$ является константой, причем константы сходятся к $e$ и к $e^2$. Получили  $e=e^2$, следовательно, $\mu\otimes\mu\,(E)\in\{0,1\}$. Таким образом,  независимый фактор тривиален. Теоремы 7.2 и 7.1 доказаны.

Из теоремы 7.2 также вытекает, что свойство системы не иметь  жестких факторов (мягкое перемешивание) является стабильным свойством.

\vspace{2mm}
 \bf Теорема 7.4. \it  Свойство мягкого перемешивания стабильно.
\rm 
\vspace{2mm} 

Доказательство. Пусть расширение $R$ мягко перемешивающего автоморфизма $S$ обладает жестким фактором,
изоморфным автоморфизму $T$,  $T^{n_i}\to I$, $n_i\to\infty$. Пусть $P$ -- 
марковский оператор, сплетающий $T$  и $S$. Имеем
$$ SP=PT, \  S^{n}Pf =P T^{n}f.$$
При $n_i\to\infty$ получаем
$$T^{n_i}f\to f, \ S^{n_i}Pf\to Pf $$
Так как $S$ не имеет жестких факторов, функция  $Pf$ обязана быть  константой. Иначе алгебра множеств, которую порождает функция $Pf$ и ее сдвиги $S^mPf$, была бы нетривиальным жестким фактором. 
Таким образом, сплетение $P$ тривиально. Это означает, что  жесткий фактор  независим относительно 
 мягко перемешивающего фактора. Мягкое перемешивание влечет за собой слабое перемешивание.
Для завершения доказательства осталось применить теорему 7.2.

\bf Кратное перемешивание и отсутствие нетривиальных джойнингов с парной независимостью. \rm
Автоморфизм $S$ перемешивает с кратностью $n$, если   для любых $A_0,A_1,\dots,A_n\in\B$ при $k_1,\dots, k_n\to\infty$ выполнено
$$\mu\left(A_0\cap S^{k_1}A_1\cap S^{k_1+k_2}A_2\dots\cap  S^{k_1+\dots+k_n}A_n\right)\to 
\mu(A_0 )\mu(A_1 )\dots \mu(A_n).$$

Самоприсоединением  порядка $n>2$ с парной независимостью называется $S\times \dots\times S$-инвариантная мера на кубе 
$X^n=X\times \dots\times X$ ($n$ сомножителей) с проекциями $\mu\otimes\mu$ на все двумерные декартовы  грани куба $X^n$.

Если для всех $n>2$  самоприсоединения с парной независимостью  для автоморфизма  $S$ тривиальны  (совпадают с $\mu^n=\mu\otimes\dots\otimes\mu$),  говорим,  что  $S$ обладает $JR$-свойством.
 Как известно,  перемешивающий автоморфизм с  $JR$-свойством обладает перемешиванием всех кратностей. С учетом установленной нами стабильности свойства перемешивания и стабильности $JR$-свойства, доказанной в  \cite{R23}, получаем следующий  факт.

\vspace{2mm}
 \bf Теорема 7.5. \it  Типичные расширения наследуют одновременное выполнение свойства кратного перемешивания и $JR$-свойства автоморфизма.
\rm 
\vspace{2mm}

\vspace{5mm} 
\bf Благодарности. \rm  Автор благодарит рецензента за замечания и выражает 
признательность  Э.~Глазнеру, Б.~Вейсу и  Ж.-П.~Тувено за полезные обсуждения.

\end{fulltext}

\end{document}